\documentclass[12pt]{amsart}
\usepackage{amssymb}

\newcommand{\Ad}{\mbox{\textnormal{Ad}}}
\newcommand{\Aut}{\mbox{\textnormal{Aut}}}
\newcommand{\A}{\mathcal{A}}
\newcommand{\AAA}{\mathfrak{A}}
\newcommand{\B}{\mathcal{B}}

\newcommand{\C}{\mathbb{C}}
\newcommand{\CC}{\mathfrak{C}}

\newcommand{\eps}{\varepsilon}
\newcommand{\f}{\mathcal{F}}
\newcommand{\F}{\mathbb{F}}
\newcommand{\h}{\mathfrak{H}}

\newcommand{\id}{\mbox{\textnormal{id}}}

\newcommand{\Inn}{\mbox{\textnormal{Inn}}}
\newcommand{\kom}{\mathcal{K}}

\newcommand{\M}{\mathcal{M}}

\newcommand{\Mt}{\mathbb{M}_2}

\newcommand{\N}{\mathcal{N}}

\newcommand{\R}{\mathbb{R}}

\newcommand{\U}{\mathcal{U}}

\newcommand{\X}{\mathfrak{X}}
\newcommand{\z}{\mathcal{Z}}
\newcommand{\Z}{\mathbb{Z}}

\newtheorem{theorem}{Theorem}
\newtheorem{proposition}[theorem]{Proposition}
\newtheorem{corollary}[theorem]{Corollary}
\newtheorem{lemma}[theorem]{Lemma}

\theoremstyle{definition}
\newtheorem{definition}[theorem]{Definition}

\newtheorem{problem}[theorem]{Problem}

\theoremstyle{remark}
\newtheorem{remark}[theorem]{Remark}

\numberwithin{equation}{section}
\numberwithin{theorem}{section}

\begin{document}

\title[Locally inner automorphisms]{Locally inner automorphisms of operator algebras}
\author{David Sherman}
\address{Department of Mathematics\\ University of Virginia\\ P.O. Box 400137\\ Charlottesville, VA 22904}
\email{dsherman@virginia.edu}
\subjclass[2000]{Primary 46L40; Secondary 47C15}
\keywords{von Neumann algebra, locally inner automorphism, diagonal sum}

\begin{abstract}
In this paper an automorphism of a unital $C^*$-algebra is said to be \textit{locally inner} if on any element it agrees with some inner automorphism.  We make a fairly complete study of local innerness in von Neumann algebras, incorporating comparison with the pointwise innerness of Haagerup-St\o rmer.  On some von Neumann algebras, including all with separable predual, a locally inner automorphism must be inner.  But a transfinitely recursive construction demonstrates that this is not true in general.  As an application, we show that the diagonal sum $(x, y) \mapsto (\begin{smallmatrix} x & 0 \\ 0 & y \end{smallmatrix})$ descends to a well-defined map on the automorphism orbits of a unital $C^*$-algebra if and only if all its automorphisms are locally inner.
\end{abstract}

\maketitle

\begin{quote}
\textit{`` `The inner truth is hidden --- luckily, luckily.  But I felt it all the same\dots .' "} (Joseph Conrad, \underline{Heart of Darkness})
\end{quote}

\section{Introduction} \label{S:intro}

In this paper we study the following concept.

\begin{definition} \label{T:def}
Let $\theta$ be an automorphism of a unital $C^*$-algebra $\AAA$.  We say that $\theta$ is \textbf{locally inner} if for every $x \in \AAA$, there is a unitary $u_x \in \AAA$ such that $\theta(x) = u_x x u_x^*$.
\end{definition}

This usage of the term ``locally inner" conflicts with some previous authors', but it is entirely in keeping with current nomenclature.  We explain this and make some basic observations in the rest of this introduction.  Section 2 of the paper is devoted to locally inner automorphisms of von Neumann algebras; Section 3 explains how locally inner automorphisms are naturally related to diagonal sums.  The last section consists of comments and open problems.

Let us first review our notation and terminology.  For a complex Hilbert space $\h$, we write $\B(\h)$ (resp. $\kom(\h)$) for the algebra of bounded (resp. compact) linear operators on $\h$.  The \textit{Calkin algebra} $\CC$ is the quotient $C^*$-algebra $\B(\ell^2)/\kom(\ell^2)$.  When not otherwise modified, $\AAA$ stands for an arbitrary unital $C^*$-algebra, and $\M$ is an arbitrary von Neumann algebra; these need not be represented on Hilbert spaces.  The algebra of $2 \times 2$ complex matrices is $\Mt$, and its matrix units are $\{e_{ij}\}_{i,j = 1,2}$.  We write $\U(\AAA)$ for the group of unitaries, $\z(\AAA)$ for the center, $\Aut(\AAA)$ for the group of automorphisms, $\Inn(\AAA)$ for the group of inner automorphisms, and $\Mt(\AAA)$ for $\Mt \otimes \AAA$.  An automorphism which is not inner will sometimes be called \textit{outer}, and ``maximal abelian *-subalgebra" is abbreviated to \textit{MASA}.  The \textit{automorphism orbit} of $x \in \AAA$ is $\{\theta(x) \mid \theta \in \Aut(\AAA)\}$.  Throughout the paper ``generate" should be understood in the von Neumann algebra sense whenever the ambient algebra is a von Neumann algebra; otherwise ``generate" has its $C^*$-meaning.  We denote the unitally included von Neumann subalgebra generated by a set or element as $W^*(\cdot)$.  The free group on $n$ generators is $\F_n$, and the von Neumann algebra generated by the left regular representation of a group $G$ is $L(G)$.  On a $\text{II}_1$ factor with tracial state $\tau$, the $L^2$ topology is generated by the norm $\|x\|_2 = \tau(x^*x)^{1/2}$, and we will use the fact that it equals the strong topology on bounded sets (see for instance \cite[Proposition III.2.2.17]{Bl}).

\smallskip

In functional analysis the adjective ``local" has many meanings, one of which is the following: a map has a property \textit{locally} if, on any input, it agrees with a map that actually has the property.  Much of the interest in this idea derives from influential work of Kadison (\cite{K1990}) and Larson-Sourour (\cite{LS}), who each showed in 1990 that certain local properties must be global.  Let us illustrate this with an explicit statement.  A \textit{local derivation} on a $C^*$-algebra $\AAA$ is a linear map $T$ from $\AAA$ to a Banach $\AAA$-bimodule $\X$ such that for every $a \in \AAA$, there is a derivation $\delta_a: \AAA \to \X$ with $T(a) = \delta_a(a)$.  Improving a theorem of Kadison, Johnson (\cite{J}) showed that local derivations on $C^*$-algebras are derivations.  If we specialize to the case where algebra and bimodule are the same von Neumann algebra $\M$, and use the fact that a derivation from $\M$ into itself must be inner (\cite{K1966,Sa}), we obtain the implication for linear maps $T: \M \to \M$,
\begin{equation} \label{E:locder}
[(\forall x)(\exists d_x)(T(x) = d_x x - x d_x)] \: \Rightarrow \: [(\exists d)(\forall x)(T(x) = dx - xd)].
\end{equation}
The reader will observe the resemblance to a statement that a locally inner automorphism of $\M$ must be inner.

As given in Definition \ref{T:def}, locally inner automorphisms seem to be a new class.  (The analogous concept in group theory has a long history.  To avoid bombarding the reader with conflicting terminology, we postpone a sketch of this history to Section \ref{S:groups}.)  But outer locally inner automorphisms have made appearances in the operator algebra literature, perhaps most notably in recent work of Phillips and Weaver (\cite{PW}).  They proved that the continuum hypothesis implies the existence of outer automorphisms of the Calkin algebra, and the automorphisms they display are locally inner by construction.  The reader will of course find other examples in this paper, the simplest in Section \ref{S:card}.

Here is a useful reformulation: an automorphism is locally inner if, on any singly-generated subalgebra, it agrees with some inner automorphism.  So when the algebra itself is singly-generated, ``locally inner" is the same as ``inner."  This actually goes a long way for von Neumann algebras.  We remind the reader of the \textit{generator problem}, which asks whether all von Neumann algebras with separable predual are singly-generated.  Although open in general (in particular, for $L(\F_3)$), it is known to have an affirmative answer for properly infinite algebras (\cite[Theorem 2]{W}) and type I algebras (\cite{Pe}).  Like the implication \eqref{E:locder}, this suggests that locally inner automorphisms of von Neumann algebras try very hard to be inner.  In the next section we study the extent to which they succeed.

Unfortunately ``locally inner" has already acquired other meanings in $C^*$-theory (\cite{L} and later variations).  Speaking in general terms, these refer to automorphisms which can be recovered by piecing together inner sub-automorphisms subordinate to a partition of unity on some topological space, typically the spectrum of the algebra.  In a sense our locally inner automorphisms can be decomposed into inner pieces by a \textit{noncommutative} partition of unity, i.e. singly-generated subalgebras.  (See the proof of Proposition \ref{T:free} for an illustration.)  However, these pieces need not leave the subalgebras invariant.

\smallskip

In Section 3 we give an application in which locally inner automorphisms arise naturally.  To motivate this, note that the usual direct sum of Hilbert space operators is ``diagonal":
$$(x_1 \in \B(\h_1), x_2 \in \B(\h_2)) \mapsto (\begin{smallmatrix} x_1 & 0 \\ 0 & x_2 \end{smallmatrix}) \in \B(\h_1 \oplus \h_2).$$
Suppose we vary $x_1$ and $x_2$ to $y_1$ and $y_2$, respectively, by automorphisms of the ambient algebras.  Since automorphisms of type I factors are inner, we must have unitaries $u_j \in \B(\h_j)$ satisfying $y_j = u_j x_j u_j^*$ for $j=1,2$, and then
\begin{equation} \label{E:bh}
(\begin{smallmatrix} y_1 & 0 \\ 0 & y_2 \end{smallmatrix}) = (\begin{smallmatrix} u_1 & 0 \\ 0 & u_2 \end{smallmatrix})(\begin{smallmatrix} x_1 & 0 \\ 0 & x_2 \end{smallmatrix})(\begin{smallmatrix} u_1 & 0 \\ 0 & u_2 \end{smallmatrix})^*.
\end{equation}
We therefore have a well-defined map on automorphism orbits,
\begin{equation} \label{E:classes}
([x_1],[x_2]) \mapsto \left[(\begin{smallmatrix} x_1 & 0 \\ 0 & x_2 \end{smallmatrix})\right].
\end{equation}

Does \eqref{E:classes} make sense in unital $C^*$-algebras other than $\B(\h)$?  We make this more precise with

\begin{definition}
The \textbf{diagonal sum} in $\AAA$ is the map
\begin{equation} \label{E:mdsum}
\AAA \times \AAA \to \Mt(\AAA), \qquad (x , y) \mapsto (\begin{smallmatrix} x & 0 \\ 0 & y \end{smallmatrix}) \overset{\mbox{\smaller[5] \textnormal{def} \larger[5]}}{=} x \oplus_d y.
\end{equation}
\end{definition}

It is important that the output belong to $\Mt(\AAA)$ and not $\AAA \oplus \AAA$.  If we use the latter, the sum of automorphism orbits is always well-defined, even when $x$ and $y$ belong to different algebras.  It turns out that \eqref{E:mdsum} induces a well-defined map on the automorphism orbits of $\AAA$ if and only if all automorphisms of $\AAA$ are locally inner.

\smallskip

This investigation was motivated by our recent work in \cite[Appendix]{S2006div}, where isomorphism classes of operators in von Neumann algebras are amplified by cardinals and, more generally, appropriate coupling functions.  The results here show that these amplifications do not typically arise as iterations of a diagonal sum.

\section{Locally inner automorphisms of von Neumann algebras} \label{S:linn}

The main goal of this section is to compare innerness and local innerness for automorphisms of von Neumann algebras.  We will see that on many von Neumann algebras, including all with separable predual, a locally inner automorphism must be inner.  But on a suitably large von Neumann algebra we are able to distinguish the two notions by a transfinitely recursive construction.

\smallskip

The predual version of local innerness has already been studied: Haagerup and St\o rmer (\cite{HSeq,HSpinn,HSinj}) defined $\theta \in \Aut(\M)$ to be \textit{pointwise inner} if for any normal state $\varphi \in \M_*^+$ there is a unitary $u_\varphi \in \M$ such that $\varphi \circ \theta = \varphi \circ \text{Ad}(u_\varphi)$.  To distinguish this property from local innerness, consider an outer modular automorphism of a type III algebra with separable predual; since the algebra is singly-generated, the automorphism is not locally inner.  On the other hand all modular automorphisms are pointwise inner (\cite[Proposition 12.6]{HSeq}).

In this paragraph we discuss a \textit{semifinite} algebra $\M$ with faithful normal semifinite trace $\tau$.  Haagerup and St\o rmer showed (\cite[Lemma 2.2]{HSpinn}) that pointwise innerness of $\theta$ is equivalent to the condition: for every $h \in \M_+ \cap L^1(\M, \tau)$, there is a unitary $u_h \in \M$ satisfying $\theta(h) = u_h h u_h^*$.  (Pointwise innerness on \textit{finite} algebras thus guarantees an implementing unitary for any self-adjoint $h$, while local innerness gives an implementing unitary for any two self-adjoints simultaneously.)  This makes it clear that
\begin{equation} \label{E:sf}
\text{inner} \quad \Rightarrow \quad \text{locally inner} \quad \Rightarrow \quad \text{pointwise inner}.
\end{equation}
Using a powerful result of Popa (\cite[Theorem 4.2]{Psing}), Haagerup and St\o rmer also proved that when $\M$ has separable predual, pointwise inner automorphisms are inner (\cite[Proposition 12.5]{HSeq}), making all conditions of \eqref{E:sf} equivalent.

Applying this result to the $\text{II}_1$ summand of an algebra with separable predual, and single-generator results to the complementary summand, we obtain

\begin{theorem} \label{T:linn}
On a von Neumann algebra with separable predual, locally inner automorphisms are inner.
\end{theorem}

Of course, if the generator problem is ever settled positively, Theorem \ref{T:linn} would become trivial.

\smallskip

We now turn to von Neumann algebras without separable predual, and we consider the converse implications to \eqref{E:sf}.  We will see that both fail, even for $\text{II}_1$ factors.

Haagerup and St\o rmer constructed outer pointwise inner automorphisms of $\text{II}_1$ factors (\cite[Theorem 6.2]{HSpinn}).  We can show that some (possibly all) of their examples fail to be locally inner.  Perhaps taxing the reader's patience for variations of ``inner," we explicitly recall that $\theta \in \Aut(\M)$ is \textit{approximately inner} if
\begin{equation} \label{E:ainn}
(\forall \eps > 0)(\forall \f \underset{\text{finite}}{\subset} \M)(\exists u \in \U(\M))(\max_{x \in \f} \|u x u^* - \theta(x)\|_2 < \eps).
\end{equation}

\begin{proposition} \label{T:pnotl}
There exists an automorphism of a $\text{II}_1$ factor which is pointwise inner but not locally inner.
\end{proposition}

\begin{proof}
We briefly review the Haagerup-St\o rmer setup and refer the reader to \cite[Theorem 6.2]{HSpinn} for full details.  Take a $\text{II}_1$ factor $\M$ with separable predual.  Choose a free ultrafilter $\omega$ on $\mathbb{N}$, and form the tracial ultraproduct $\M^\omega$, which is also a $\text{II}_1$ factor.  Finally take $\theta \in \Aut(\M)$, and let $\theta^\omega \in \Aut(\M^\omega)$ be the induced automorphism.  Haagerup and St\o rmer showed (1) $\theta^\omega$ is pointwise inner; (2) if $\theta^\omega$ is inner, then $\theta$ must be approximately inner.  Choosing non-approximately inner $\theta$ then gives outer pointwise inner $\theta^\omega$.

Next we show that (3) if $\M$ has a generator and $\theta^\omega$ is locally inner, then $\theta$ must be approximately inner.  Thus in the Haagerup-St\o rmer examples of the previous paragraph, $\theta^\omega$ is not locally inner, either -- at least if we know that $\M$ is singly-generated.  And sometimes we do: for example, take $\M$ to be $L(\F_2)$ and $\theta$ to be the non-approximately inner automorphism induced by permuting the group generators (\cite[Exercice III.7.15]{Dbook}).

Here is the proof of (3).  Let $y$ be a generator for $\M$, and let $\pi$ be the embedding of $\M$ in $\M^\omega$ as (cosets of) constant sequences.  The local innerness of $\theta^\omega$, applied to $\pi(y)$, says that there is a unitary $U \in \M^\omega$ such that $\theta^\omega(\pi(y)) = U \pi(y) U^*$.  Since $U$ has a representing sequence of unitaries $(u_n)$, we have that
\begin{equation} \label{E:conv}
\|u_n y u_n^* - \theta(y)\|_2 \to 0 \qquad \text{as }n \to \omega.
\end{equation}
Clearly we can replace $y$ in \eqref{E:conv} by any polynomial in $y$ and $y^*$.  By the strong density of these polynomials in $W^*(y) = \M$ and use of the triangle inequality, it then follows that \eqref{E:conv} is true for all elements of $\M$.  Now suppose $\eps > 0$ and a finite set $\f \subset \M$ are given.  For each $x \in \f$ there is an open neighborhood $G_x$ of $\omega$ in the Stone-\v{C}ech compactification $\beta \mathbb{N}$ such that $n \in G_x$ implies $\|u_n x u_n^* - \theta(x)\|_2 < \eps$.  The finite intersection of the $G_x$ is an open neighborhood of $\omega$, in particular non-empty, and so contains an integer $m$.  Evidently $\max_{x \in \f} \|u_m x u_m^* - \theta(x)\|_2 < \eps$ as desired.
\end{proof}

\begin{remark}
Here is the underlying point of the preceding paragraph.  Any non-approximately inner automorphism has a witness $(\eps, \f)$ violating \eqref{E:ainn}.  When $\M$ is singly-generated, we can always take $\f$ to consist of a single element (the generator).  Are there any non-approximately inner automorphisms for which $\f$ must consist of at least two elements?  This would be required to distinguish ``locally inner" from ``inner" by an example of Haagerup-St\o rmer type.

If the generator problem has a positive solution, then any non-approximately inner automorphism has a witness $(\eps,\f)$ in which $\f$ is a singleton, even when $\M$ is not assumed to have separable predual.  To see this, suppose $(\eps,\f)$ is a witness for $\theta$.  Set
$$\N = W^*(\{\theta^n(x) \mid n \in \Z, x \in \f\}) \subseteq \M.$$
Note that $\theta$ leaves $\N$ invariant, and its restriction is apparently not approximately inner.  But $\N$ is countably generated, so by our assumption, it has a generator.  This generator must witness the non-approximate innerness of $\theta$.
\end{remark}

To show that local innerness implies innerness in a specific algebra, one has to get some global control on an automorphism through its action on singly-generated subalgebras.  This requires a certain type of ``rigid" inclusion, and a prime example is the Popa result behind Theorem \ref{T:linn}.  There one finds a singly-generated MASA $\A$ such that the automorphism is determined, up to conjugation by a unitary from $\A$, by its action on $\A$.  One obstacle to a general application of this technique is that some large algebras like $\M^\omega$ have no singly-generated MASAs at all (\cite[Proposition 4.3]{Pkad}).  On the other hand, arbitrarily large free group factors have singly-generated MASAs with enough rigidity for the following

\begin{proposition} \label{T:free}
On a free group factor with any (possibly uncountable) number of generators, a locally inner automorphism must be inner.
\end{proposition}

\begin{proof}
Let $\{g_j\}_{j \in J}$ be the generators of the free group $\F_J$, where $J$ is an arbitrary index set.  Writing $\lambda_{g_j} \in L(\F_J)$ for the unitary on $\ell^2(\F_J)$ which left-convolves with $g_j$, set $h_j = i \, \text{Log} (\lambda_{g_j})$.  So $L(\F_J)$ is generated by the self-adjoint operators $\{h_j\}$, and the subalgebras $\A_j = W^*(h_j)$ are abelian, singular (\cite[Remark 6.3]{Porth}), and free from each other.  (Recall that an abelian subalgebra is \textit{singular} if it contains its normalizer.)

Supposing that $\theta \in \Aut(L(\F_J))$ is locally inner, we condense notation somewhat and write
$$\theta(h_j + i h_k) = u_{jk} (h_j + i h_k) u_{jk}^*, \qquad j,k \in J.$$
Choose any three indices $j,k,l$.  By singularity,
$$u_{jk} h_j u_{jk}^* = \theta(h_j) = u_{jl} h_j u_{jl}^* \Rightarrow \text{Ad}(u_{jk}^* u_{jl})(h_j) = h_j \Rightarrow u_{jk}^* u_{jl} \in \A_j.$$
Consider the equation
\begin{equation} \label{E:key}
u_{jk}^* u_{kl} = (u_{jk}^* u_{jl})(u_{jl}^* u_{kl}).
\end{equation}
The left-hand side of \eqref{E:key} belongs to $\A_k$, while the right-hand side belongs to $\A_j \cdot \A_l$.  By freeness $\A_k \cap (\A_j \cdot \A_l) = \C$, so the quantity in \eqref{E:key} is a scalar.   Looking again at the left-hand side, we now deduce $\text{Ad}(u_{jk}) = \text{Ad}(u_{kl})$.  So for arbitrary $l$, we have shown
$$\theta(h_l) = \text{Ad}(u_{kl})(h_l) = \text{Ad}(u_{jk})(h_l).$$
We conclude that $\theta = \text{Ad}(u_{jk})$.
\end{proof}

The proof above is ``cohomological": the algebras $W^*(h_j + i h_k)$ play the role of an open cover, with the unitaries $u_{jk}^* u_{jl}$ as 1-cocycles.  The argument has much in common with \cite[Section 4]{L}.

The following construction was inspired by the Akemann-Weaver counterexample to Naimark's problem (\cite{AW}).

\begin{theorem} \label{T:ctrex}
There exists an automorphism of a $\text{II}_1$ factor which is locally inner but not inner.
\end{theorem}

\begin{proof}
Our construction is transfinitely recursive.  To each ordinal $\alpha < \aleph_1$ we will associate a hyperfinite $\text{II}_1$ factor $\M_\alpha$ with separable predual and distinguished self-adjoint unitary $u_\alpha$.  When $\beta < \alpha$ we will require that
\begin{enumerate}
\item $\M_\beta \subsetneq \M_\alpha$ (unital inclusion),
\item $\text{Ad}(u_\alpha)\mid_{\M_\beta} = \text{Ad}(u_\beta)\mid_{\M_\beta}$, and
\item there is no $u \in \U(\M_\beta)$ satisfying $\text{Ad}(u) \mid_{\M_\alpha}= \text{Ad}(u_\alpha) \mid_{\M_\alpha}$.
\end{enumerate}
Assuming momentarily that this has been done, we show how to complete the proof.

We first claim that the union $\M = \cup_{\alpha < \aleph_1} \M_\alpha$ is already a von Neumann algebra.  Note that there is a unique tracial state on each $\M_\alpha$, so these must be coherent.  The union of these traces defines a tracial state $\tau$ on $\M$, and the strong topology on the unit ball in the faithful GNS representation for $(\M, \tau)$ is determined by the $L^2$ norm.  By the Kaplansky density theorem, $\M$ is a von Neumann algebra if and only if its norm-closed unit ball $\M$ is strongly complete, so it suffices to check $L^2$ completeness.  This allows us, significantly, to work with sequences instead of arbitrary nets.  So let $\{x_n\} \subset \M$ be contractions which are Cauchy in $L^2$.  Each $x_n$ belongs to some $\M_{\alpha_n}$, and $\sup \alpha_n < \aleph_1$.  It follows that the $x_n$ all belong to $\M_{\sup \alpha_n} \subset \M$ and have a strong limit there.  This shows that $\M$ is a von Neumann algebra -- in fact it is a $\text{II}_1$ factor, since it admits a unique tracial state.

Let $\theta$ be the automorphism of $\M$ such that $\theta(x) = u_\alpha x u_\alpha^*$ for $x \in \M_\alpha$.  This is well-defined by (2), and evidently it is locally inner.  But it is not inner, as any unitary belongs to some $\M_\alpha$ and then cannot implement $\theta$ on $\M_{\alpha + 1}$ by (3).

We now construct the $(\M_\alpha, u_\alpha)$.  Let $u_0$ be any self-adjoint unitary in $\M_0$, a hyperfinite $\text{II}_1$ factor with separable predual.  If $\alpha$ is a successor ordinal, say $\alpha = \beta + 1$, set
$$\M_\alpha = \Mt \otimes \M_\beta \supsetneq \mathbb{C} \otimes \M_\beta = \M_\beta; \quad u_\alpha = (\begin{smallmatrix} 0 & 1 \\ 1 & 0 \end{smallmatrix}) \otimes u_\beta.$$
If $\alpha$ is a limit ordinal, first define $\N_\alpha$ to be the von Neumann algebra direct limit of $\{\M_\beta \mid \beta < \alpha\}$ with respect to the coherent trace, i.e. the strong closure in the GNS representation for $\cup_{\beta < \alpha} \M_\beta$ with respect to this trace.  The automorphisms $\text{Ad}(u_\beta) \mid_{\M_\beta}$ agree on their common domains and so define a $\|\cdot\|_2$-continuous automorphism on $\cup_{\beta < \alpha} \M_\beta$ which extends to a period two automorphism $\sigma_\alpha$ of $\N_\alpha$.  We enlarge the algebra so that $\sigma_\alpha$ can be implemented by a unitary.  Set $\M_\alpha = \Mt(\N_\alpha)$, and choose the inclusion and unitary as follows:
$$\N_\alpha \ni x \mapsto (\begin{smallmatrix} x & 0 \\ 0 & \sigma_\alpha(x) \end{smallmatrix}) \in \M_\alpha, \qquad u_\alpha = (\begin{smallmatrix} 0 & 1_{\N_\alpha} \\ 1_{\N_\alpha} & 0 \end{smallmatrix}).$$
It is straightforward to verify that the $(\M_\alpha, u_\alpha)$ have all the desired properties.
\end{proof}

\begin{remark}
For limit ordinals $\alpha$, the last construction is nothing other than
$$\M_\alpha = (\N_\alpha \rtimes_{\sigma_\alpha} \Z_2) \rtimes_{\hat{\sigma}_\alpha} \Z_2,$$
the isomorphism with $\Mt(\N_\alpha) \cong \N_\alpha \otimes \B(\ell^2(\Z_2))$ being a specific case of Takesaki duality (\cite[Theorem 4.5]{T1973}).
\end{remark}

\section{Diagonal sums} \label{S:ds}

Widening our scope back to $C^*$-algebras, here we consider whether the diagonal sum \eqref{E:mdsum} makes sense for automorphism orbits.  For convenience some versions of the desired property are collected in an easy lemma.  (To prove (3) $\to$ (2), apply $\text{id} \otimes \sigma^{-1}$ to $(\begin{smallmatrix} \theta(x) & 0 \\ 0 & \sigma(y) \end{smallmatrix})$ and rename $\sigma^{-1} \circ \theta$ as $\theta$.)

\begin{lemma} \label{T:prob}
For a unital $C^*$-algebra $\AAA$, the following properties are equivalent:
\begin{enumerate}
\item there is a well-defined diagonal sum operation on automorphism orbits of operators in $\AAA$,
$$([x],[y]) \mapsto [x \oplus_d y] \overset{\mbox{\smaller[5] \textnormal{def} \larger[5]}}{=} [x] \oplus_d [y];$$
\item for any $\theta, \sigma \in \Aut(\AAA)$ and $x,y \in \AAA$, there is $\gamma \in \Aut(\Mt(\AAA))$ with
$$\gamma(\begin{smallmatrix} x & 0 \\ 0 & y \end{smallmatrix}) = (\begin{smallmatrix} \theta(x) & 0 \\ 0 & \sigma(y) \end{smallmatrix});$$
\item for any $\theta \in \Aut(\AAA)$ and $x,y \in \AAA$, there is $\gamma \in \Aut(\Mt(\AAA))$ with
$$\gamma(\begin{smallmatrix} x & 0 \\ 0 & y \end{smallmatrix}) = (\begin{smallmatrix} \theta(x) & 0 \\ 0 & y \end{smallmatrix}).$$
\end{enumerate}
\end{lemma}

The following lemma is straightforward.

\begin{lemma} \label{T:inv}
Let $\gamma$ be an automorphism of $\Mt(\AAA)$ which leaves $(\begin{smallmatrix} 0 & 0 \\ 0 & 1 \end{smallmatrix})$ invariant.  Then there are $\sigma \in \Aut(\AAA)$ and $u \in \U(\AAA)$ such that
\begin{equation} \label{E:factor}
\gamma = \Ad(\begin{smallmatrix} u & 0 \\ 0 & 1 \end{smallmatrix}) \circ (\id \otimes \sigma).
\end{equation}
\end{lemma}

\begin{proof}
Since $\gamma$ restricts to an automorphism of $(e_{22} \otimes 1)\Mt(\AAA)(e_{22} \otimes 1) \cong \AAA$, we may define $\sigma$ by
$$\gamma(\begin{smallmatrix} 0 & 0 \\ 0 & x \end{smallmatrix}) = (\begin{smallmatrix} 0 & 0 \\ 0 & \sigma(x) \end{smallmatrix}), \qquad \forall x \in \AAA.$$

Notice
$$\gamma(\begin{smallmatrix} 0 & 1 \\ 0 & 0 \end{smallmatrix}) \gamma(\begin{smallmatrix} 0 & 1 \\ 0 & 0 \end{smallmatrix})^* = \gamma(\begin{smallmatrix} 1 & 0 \\ 0 & 0 \end{smallmatrix}) = (\begin{smallmatrix} 1 & 0 \\ 0 & 0 \end{smallmatrix}), \quad \gamma(\begin{smallmatrix} 0 & 1 \\ 0 & 0 \end{smallmatrix})^* \gamma(\begin{smallmatrix} 0 & 1 \\ 0 & 0 \end{smallmatrix}) = \gamma(\begin{smallmatrix} 0 & 0 \\ 0 & 1 \end{smallmatrix}) = (\begin{smallmatrix} 0 & 0 \\ 0 & 1 \end{smallmatrix}).$$
Thus $\gamma(\begin{smallmatrix} 0 & 1 \\ 0 & 0 \end{smallmatrix})$ is a partial isometry between $(\begin{smallmatrix} 1 & 0 \\ 0 & 0 \end{smallmatrix})$ and $(\begin{smallmatrix} 0 & 0 \\ 0 & 1 \end{smallmatrix})$, necessarily of the form $(\begin{smallmatrix} 0 & u \\ 0 & 0 \end{smallmatrix})$ for some unitary $u \in \AAA$.

It suffices to verify \eqref{E:factor} on any element of the form $e_{ij} \otimes x$ ($x \in \AAA$, $i,j \in \{1,2\}$), and this is straightforward algebra.  For example,
\begin{align*}
\gamma(e_{11} \otimes x) &= \gamma[(e_{12} \otimes 1)(e_{22} \otimes x)(e_{21} \otimes 1)] \\ &= (e_{12} \otimes u)(e_{22} \otimes \sigma(x))(e_{21} \otimes u^*) \\ &= e_{11} \otimes u \sigma(x) u^* \\ &= [\text{Ad}(\begin{smallmatrix} u & 0 \\ 0 & 1 \end{smallmatrix}) \circ (\text{id} \otimes \sigma)] (e_{11} \otimes x). \qedhere
\end{align*}
\end{proof}

We now analyze the consequences of the condition in Lemma \ref{T:prob}(3) for an individual automorphism $\theta$.

\begin{theorem} \label{T:main}
For an automorphism $\theta$ of a unital $C^*$-algebra $\AAA$, the following conditions are equivalent:
\begin{enumerate}
\item for any $x,y \in \AAA$, there is $\gamma \in \Aut (\Mt (\AAA))$ such that
$$\gamma(\begin{smallmatrix} x & 0 \\ 0 & y \end{smallmatrix}) = (\begin{smallmatrix} \theta(x) & 0 \\ 0 & y \end{smallmatrix});$$
\item $\theta$ is locally inner.
\end{enumerate}
\end{theorem}

\begin{proof}
The implication (2) $\to$ (1) is easy: set $\gamma = \text{Ad}(\begin{smallmatrix} u_x & 0 \\ 0 & 1 \end{smallmatrix})$.

Now suppose that $\theta$ satisfies (1).  Choose $x \in \AAA$; after rescaling we may assume $\|x\|<1$.  By hypothesis there is $\gamma \in \Aut(\AAA)$ such that $\gamma(\begin{smallmatrix} x & 0 \\ 0 & x+2 \end{smallmatrix}) = (\begin{smallmatrix} \theta(x) & 0 \\ 0 & x+2 \end{smallmatrix})$.

The spectrum of $(\begin{smallmatrix} x & 0 \\ 0 & x+2 \end{smallmatrix})$ is the disjoint union of the spectra of $x$ and $x+2$, which are contained in open unit disks centered at 0 and 2 respectively.  Let $f$ be a function, analytic on a neighborhood of the spectrum of $(\begin{smallmatrix} x & 0 \\ 0 & x+2 \end{smallmatrix})$, which is 0 on the spectrum of $x$ and 1 on the spectrum of $x+2$.  By Riesz holomorphic functional calculus,
$$\gamma(\begin{smallmatrix} 0 & 0 \\ 0 & 1 \end{smallmatrix}) = \gamma(f(\begin{smallmatrix} x & 0 \\ 0 & x+2 \end{smallmatrix})) =f(\gamma(\begin{smallmatrix} x & 0 \\ 0 & x+2 \end{smallmatrix})) = f(\begin{smallmatrix} \theta(x) & 0 \\ 0 & x+2 \end{smallmatrix}) = (\begin{smallmatrix} 0 & 0 \\ 0 & 1 \end{smallmatrix}).$$
Lemma \ref{T:inv} locates $\sigma \in \Aut(\AAA)$ and $u \in \U(\AAA)$ such that $\gamma = \Ad(\begin{smallmatrix} u & 0 \\ 0 & 1 \end{smallmatrix}) \circ (\id \otimes \sigma)$.  We compute
\begin{equation} \label{E:agree}
(\begin{smallmatrix} \theta(x) & 0 \\ 0 & x+2 \end{smallmatrix}) = \gamma(\begin{smallmatrix} x & 0 \\ 0 & x+2 \end{smallmatrix}) = (\begin{smallmatrix} u\sigma(x)u^* & 0 \\ 0 & \sigma(x+2) \end{smallmatrix}).
\end{equation}
Equating the lower-right corners, we conclude that $\sigma(x)=x$.  Equating the upper-left corners then gives $\theta(x) = uxu^*$ as required.
\end{proof}

\begin{remark} \label{T:global} The proof shows that if the map $\gamma$ of Theorem \ref{T:main}(1) is not allowed to vary with $x$ and $y$, then $\theta$ must actually be inner.
\end{remark}

\begin{corollary} \label{T:diag}
A unital $C^*$-algebra $\AAA$ admits a diagonal sum on automorphism orbits (i.e., satisfies the conditions in Lemma \ref{T:prob}) if and only if all automorphisms of $\AAA$ are locally inner.
\end{corollary}

Pairing Theorem \ref{T:linn} with Corollary \ref{T:diag} produces

\begin{corollary} \label{T:autinn}
Let $\M$ be a von Neumann algebra with separable predual.  Then $\M$ admits a diagonal sum on automorphism orbits if and only if all automorphisms of $\M$ are inner.
\end{corollary}

\begin{corollary} For von Neumann algebras with separable predual and well-defined type, we can say the following.
\begin{enumerate}
\item A type I algebra admits a diagonal sum on automorphism orbits if and only if it is a direct product of nonisomorphic factors.
\item Some type II algebras admit a diagonal sum on automorphism orbits, and some do not.
\item No type III algebra admits a diagonal sum on automorphism orbits.
\end{enumerate}
\end{corollary}

\begin{proof}
(1) An automorphism fixes each type $\text{I}_k$ summand globally.  Each of these has the form $L^\infty(X_k,\mu_k) \overline{\otimes} \B(\h_k)$ ($\dim \h_k = k$) and possesses an outer automorphism exactly when $L^\infty(X_k,\mu_k) \ncong \C$.

(2) Any von Neumann algebra of the form $\N \overline{\otimes} \A$, where $\N$ is type II and $\A$ is nontrivial abelian, admits an outer automorphism of the form $\text{id} \otimes \sigma$.  Moreover, outer automorphisms on $\text{II}_1$ factors have been known since the 1950s (\cite[Remarque 1]{D}; these first examples include the generator swap on $L(\F_2)$ and specially constructed automorphisms of the hyperfinite $\text{II}_1$ factor).  The first type II algebra with \textit{no} outer automorphisms was recently displayed in \cite{IPP}.

(3) Let $\varphi$ be a faithful normal state on a type III algebra.  For some $t \in \R$, the modular automorphism $\sigma_t^\varphi$ is outer (\cite[Th\'{e}or\`{e}me 1.3.4(b)]{C1973}).
\end{proof}

What is the scope of Corollary \ref{T:diag}?  We have just discussed von Neumann algebras which lack outer automorphisms; in $C^*$-algebras there are also separable infinite-dimensional examples.  (The only simple ones are of the form $\kom(\h)$, but there are others -- see \cite[Introduction]{PW}.)  One wants to know if local innerness really adds anything to the picture.

\begin{proposition} \label{T:ex}
There is a unital $C^*$-algebra on which all automorphisms are locally inner, but not all automorphisms are inner.
\end{proposition}

\begin{proof}
Let $\{\h_j\}$ be an uncountable set of nonisomorphic Hilbert spaces.  Let $\AAA \subset \Pi_{j \in J} \B(\h_j)$ be the $C^*$-subalgebra which consists of sequences equaling some constant multiple of the identity off a countable set.  (This is bigger than the unitization of the direct sum $\Sigma_{j \in J} \B(\h_j)$, since we do not require any decay.)

Any $\alpha \in \Aut(\AAA)$ permutes the minimal central projections.  Since the $\B(\h_j)$ summands are mutually nonisomorphic, this permutation must be trivial, so that $\alpha$ fixes each $\B(\h_j)$ globally.  Therefore there are $u_j \in \U(\B(\h_j))$ such that $\alpha = \text{Ad}((u_j))$.  (Here $(u_j)$ may be viewed as an element of the multiplier algebra $\Pi_{j \in J} \B(\h_j)$.)  Now $\alpha$ is only inner if countably many $u_j$ are non-scalar, and this need not be the case.  But $\alpha$ must be locally inner: for any $(x_j) \in \AAA$,
$$\alpha(x_j) = \text{Ad}((v_j))(x_j), \quad \text{where }(v_j) \in \AAA, \quad v_j = \begin{cases} u_j, & x_j \notin \C \\ 1, & \text{otherwise}. \end{cases} \qedhere$$
\end{proof}

\section{Problems and comments}

For von Neumann algebras our understanding of locally inner automorphisms is reasonably sharp, and the remaining gaps are mostly related to the generator problem.  But for $C^*$-algebras the situation is less satisfactory.  In this final section we discuss two aspects of our examples which might be improved, then briefly sketch the analogous results from group theory, and conclude with some comments on orbits and the connection to extension theory.

\subsection{Minimizing cardinality} \label{S:card}

The reader has probably noticed that all our examples of outer locally inner automorphisms act on very large algebras.  For von Neumann algebras we have shown that this is unavoidable (Theorem \ref{T:linn}), and the construction in Theorem \ref{T:ctrex} is no bigger than it has to be.  But for $C^*$-algebras the necessity is not clear.

\begin{problem} \label{T:sep} ${}$
\begin{enumerate}
\item Find a ``minimal" $C^*$-algebra satisfying Proposition \ref{T:ex}.  Here ``minimal" could refer to the cardinality of a dense subset or to the dimension of a Hilbert space on which the algebra is faithfully represented.
\item Find a ``minimal" $C^*$-algebra which admits an outer, locally inner automorphism.
\end{enumerate}
\end{problem}

To provide more context for Problem \ref{T:sep}, we give some discussion of the exact sizes of our examples.  Readers who have no interest in distinguishing uncountable cardinals may guiltlessly skip the next paragraph.

Recall the Phillips-Weaver result: assuming the continuum hypothesis, the Calkin algebra $\CC$ admits an outer locally inner automorphism (\cite{PW}).  (Whether \textit{all} automorphisms of $\CC$ would be locally inner is left open.  Farah (\cite{Fa}) showed that in the absence of the continuum hypothesis, it is consistent with ZFC that all automorphisms of $\CC$ are inner.)  Now $\mathfrak{c} = 2^{\aleph_0}$ is the minimal cardinality of a dense set in $\CC$, as well as the minimal dimension of a Hilbert space on which $\CC$ is faithfully represented.\footnote{The size of a minimal dense set is readily deduced from the fact that $\CC$ both contains $\B(\ell^2)$ and is a quotient of it.  Concerning dimensions of representation spaces, the lower bound of $\mathfrak{c}$ results from the existence of $\mathfrak{c}$ nonzero pairwise orthogonal projections in $\CC$, and a faithful representation with dimension $\mathfrak{c}$ was constructed in the original paper by Calkin (\cite{Ca}).  Calkin worked ``with bare hands," as his paper predates the Gelfand-Naimark article \cite{GN} launching abstract $C^*$-algebra theory.  In modern language one observes that the GNS construction for any singular state on $\B(\ell^2)$ gives rise to a faithful representation of $\CC$ on a space of dimension $\mathfrak{c}$.}  Even without assuming the continuum hypothesis, one can construct an outer locally inner automorphism on a $C^*$-algebra which is no bigger, as follows.  Imitate the construction of Proposition \ref{T:ex}, but let $\{\h_j\}$ be a collection of $\aleph_1$ Hilbert spaces of dimension 2.  The resulting algebra $\AAA$ admits an outer locally inner automorphism, for example $\Ad ((u_j))$ with all $u_j = (\begin{smallmatrix} 0 & 1 \\ 1 & 0 \end{smallmatrix})$.  (It does not satisfy the condition of Proposition \ref{T:ex}, since there are non-locally inner automorphisms which permute the summands.)  One can check that $\mathfrak{c}$ is the minimal cardinality of a dense set,\footnote{Make the following choices: a countably infinite set of indices, an element of a dense set in $\Mt$ for each of these indices, and a single element of a dense set of $\C$ for all the remaining indices.  This can be done in $\aleph_1^{\aleph_0} \times \aleph_0^{\aleph_0} \times \aleph_0 = \mathfrak{c}$ ways; remove repetitions and unbounded sequences to get a dense set in $\AAA$.  But $\AAA$ also contains $\aleph_1^{\aleph_0} = \mathfrak{c}$ projections in which $\aleph_0$ entries are the identity and the rest are zero.   Since any two such projections are at distance one from each other, $\mathfrak{c}$ is also a lower bound for the cardinality of a dense set.  (For basic facts on the arithmetic and combinatorics of infinite cardinals, the reader may see \cite[Chapter 5]{Je}.)} and evidently $\AAA$ acts faithfully on a Hilbert space of dimension $\aleph_1 \times 2 = \aleph_1$.

Modulo the continuum hypothesis, then, Problem \ref{T:sep}(2) can be piquantly restated:

\begin{quote}
\textit{Must a locally inner automorphism of a separable $C^*$-algebra be inner?}
\end{quote}

\subsection{$n$-local innerness} \label{S:nloc}

All our examples of locally inner automorphisms in fact have a stronger property: on any \textit{countable} set of inputs, they agree with some inner automorphism.  We call this \textit{$\aleph_0$-local innerness}, and generally write ``$n$-local innerness" for the analogous notion corresponding to any cardinal $n$.  This follows usage from \cite{Se}, where a map is said to have a property \textit{2-locally} if, on any two inputs, it agrees with a map that actually has the property.  We have not seen higher cardinals in the literature, presumably because many basic properties -- linearity, multiplicativity, the derivation identity, etc. -- are 3-local by definition.

For cardinals $m < n$, can we distinguish $m$-local innerness from $n$-local innerness?  If $n > \aleph_0$, this can be done with an automorphism of a von Neumann algebra.  For $m$ finite, use Theorem \ref{T:ctrex}; for $m = \aleph_\alpha$, just alter the construction of Theorem \ref{T:ctrex} by replacing $\aleph_1$ with $\aleph_{\alpha + 1}$.  (This relies on the regularity of the cardinal $\aleph_{\alpha + 1}$ (\cite[Corollary 5.3]{Je}).)  But if $n \leq \aleph_0$, it seems difficult to decide whether the terms are equivalent for von Neumann algebras.  In particular, their disagreement would entail a negative solution to the generator problem.  The issue may be more tractable in $C^*$-algebras, where ``$m$-generated" and ``$n$-generated" are different even for abelian algebras.

\begin{problem}
Let $m < n \leq \aleph_0$.  Is there an $m$-locally inner automorphism which is not $n$-locally inner?
\end{problem}

\subsection{Analogues in group theory} \label{S:groups}

For several of the main ideas in the present article, the group theoretic analogue is ground run over long ago.  Unfortunately the terminology is somewhat different from ours and not always consistent from paper to paper.  As mentioned, here we have opted to follow the nomenclature from recent work in operator algebras.  It seems worthwhile to give the main terms and a brief synopsis from the group theoretic side, focusing on results which resemble the problems posed in the previous two subsections.

A group automorphism is (most commonly) called \textit{class-preserving} if it has the property that on any element it agrees with some inner automorphism.  Already Burnside's 1911 book (\cite[Note B]{B1}) asked if a class-preserving automorphism of a finite group must be inner; he answered his own question with a counterexample in 1913 (\cite{B2}).  By now several general methods have been developed for constructing counterexamples, and the order of the group has gotten as small as 32 (\cite{Wa}).  Interestingly, there is no counterexample in a finite simple group -- this was shown by referring to the classification theorem and exhausting all possibilities (\cite[Theorem C]{FS}).  %Burnside had also claimed that the quotient of the class-preserving automorphisms by the inner automorphisms must be abelian, which is not true even for finite groups.  It is at least solvable for finite (but not infinite) groups (\cite{Sah}).

A group automorphism is (most commonly) called \textit{$n$-inner} if on any set of fewer than $n$ elements, it agrees with some inner automorphism.  For any $n$, there is a group automorphism which is $n$-inner but not $(n+1)$-inner, and the group may be taken finite when $n$ is finite (\cite{N}).  In older literature (\cite{G}) $\aleph_0$-inner group automorphisms have been called \textit{locally inner}.

\subsection{Orbits and extensions}

%\smallskip

%\textbf{1.} The main results of Section \ref{S:ds} apply to nonunital $C^*$-algebras as well: just revise Definition \ref{T:def} so that $u_x$ belongs to the multiplier algebra of $\AAA$.

%\textbf{2.} Several of our conclusions (Remark \ref{T:global}, Theorem \ref{T:main} and its pairing with Theorem \ref{T:linn}) bear a pleasing superficial resemblance to a well-known result of Connes (\cite[Theorem 3.1]{C1976}): an automorphism $\alpha$ of a $\text{II}_1$ factor $\M$ with separable predual and trace $\tau$ is approximately inner if and only if there is an automorphism of the $C^*$-algebra generated by $\M \cup \M' \subseteq \B(L^2(\M, \tau))$, implementing $\alpha$ on $\M$ and the identity on $\M'$.

For any subgroup $\Pi < \Aut(\AAA)$, one can define the $\Pi$-orbit of $x \in \AAA$ as $\{\pi(x) \mid \pi \in \Pi\}$.  If we also choose $\Pi' < \Aut(\Mt(\AAA))$, we can ask whether the diagonal sum of $\Pi$-orbits is always well-defined as a $\Pi'$-orbit.

Theorem \ref{T:main} tells us that the group of locally inner automorphisms is the largest allowable choice for $\Pi$, and just as in \eqref{E:bh}, this choice is acceptable as long as $\Pi'$ contains all automorphisms of the form $\Ad (\begin{smallmatrix} u_1 & 0 \\ 0 & u_2 \end{smallmatrix})$.  But the orbit of $x \in \AAA$ under locally inner automorphisms is no different than its unitary orbit $\{uxu^* \mid u \in \U(\AAA)\}$.  So in any case, \textit{unitary orbits are the largest orbits which admit a diagonal sum.}

\smallskip

Diagonal sums will be familiar to many readers from the theory of extensions, pioneered by Busby (\cite{Bus}) and Brown-Douglas-Fillmore (\cite{BDF1,BDF2}).  In fact orbit considerations are in play there as well.  We do not review the general setup here, specializing instead to the BDF theory covering extensions of $C(X)$ by $\kom(\ell^2)$.

The main objects, called \textit{Busby invariants}, are *-homomorphisms from $C(X)$ to $\CC$.  One naturally considers equivalence classes which are orbits obtained by postcomposing with all elements from a group $\Pi < \Aut(\CC)$.   It is important that the diagonal sum of *-homomorphisms descend to a well-defined binary operation on $\Pi$-orbits; this creates a semigroup structure on classes of extensions (using $\Mt(\CC) \simeq \CC$).  For ``weak equivalence" one takes $\Pi$ to be the inner automorphisms of $\CC$, while for ``strong equivalence" $\Pi$ is the subgroup in which implementing unitaries are cosets of unitaries from $\B(\ell^2)$ .  In \cite[Remark 1.6(ii)]{BDF2}, the authors do raise other possibilities for $\Pi$, and they mention putative outer automorphisms of $\CC$ which agree with an inner automorphism on any separable commutative subalgebra.  Of course including them in $\Pi$ would give the same ``weak equivalence" relation.  (Phillips and Weaver explicitly point out that their outer automorphisms of $\CC$, being $\aleph_0$-locally inner, have no consequences for extensions (\cite[Introduction]{PW}).)  Fillmore also considers other choices for $\Pi$ in his development of the $\text{II}_\infty$ theory of extensions (\cite[Remark 1.3]{F}), but he rejects them because ``it is not clear that addition is well-defined."

It would be interesting to know if unitary orbits of extensions are the largest orbits which admit a diagonal sum.  This would mean that an appropriate version of local innerness (agreeing with an inner automorphism on the range of any *-homomorphism) is required from the members of $\Pi$.  %Actually this question makes sense in an even broader context, but its development is postponed until the answers are better understood.

%\textbf{3.} One can also ask about diagonal sums of \textit{closures} of orbits.  This is much different, and there are many positive results which do not require separable predual or an absence of outer automorphisms.  We will only make one explicit statement: there is a well-defined diagonal sum on norm closures of unitary orbits of Hilbert space operators.  This follows fairly directly, for example, from Hadwin's variation on Voiculescu's theorem characterizing approximate equivalence in terms of rank (\cite[Theorem 3.14]{H1981}).  Various other statements for $\B(\h)$ can be derived from Hadwin's work; analogues for normal operators in other von Neumann algebras follow from \cite{S2005}.

\smallskip

\textbf{Acknowledgments.} Thanks are due to Erling St\o rmer, Leonard Scott, Nik Weaver, and Ron Douglas for specific suggestions, and especially to Chuck Akemann for many helpful conversations.  This paper is dedicated to the memory of Brother John Wills (1926-2006), a mathematics teacher for the ages.


\begin{thebibliography}{BDF2}

\bibitem[AW]{AW}
C. Akemann and N. Weaver,
\emph{Consistency of a counterexample to Naimark's problem},
Proc. Natl. Acad. Sci. USA \textbf{101} (2004), 7522--7525.

\bibitem[Bl]{Bl}
B. Blackadar,
Operator Algebras,
Springer-Verlag, Berlin, 2006.

\bibitem[BDF1]{BDF1}
L. G. Brown, R. Douglas, and P. A. Fillmore,
\emph{Unitary equivalence modulo the compact operators and extensions of $C^*$-algebras}, in: Proceedings of a Conference on Operator Theory (Dalhousie Univ., Halifax, 1973), pp. 58--128, Lecture Notes in Math., vol. 345, Springer, Berlin, 1973.

\bibitem[BDF2]{BDF2}
L. G. Brown, R. Douglas, and P. A. Fillmore,
\emph{Extensions of $C^*$-algebras and $K$-homology},
Ann. of Math. (2) \textbf{105} (1977), 265--324.

\bibitem[B1]{B1}
W. Burnside,
Theory of Groups of Finite Order,
Cambridge University Press, Cambridge, 1911.

\bibitem[B2]{B2}
W. Burnside,
\emph{On the outer automorphisms of a group},
Proc. London Math. Soc. (2) \textbf{11} (1913), 40--42.

\bibitem[Bus]{Bus}
R. C. Busby,
\emph{Double centralizers and extensions of C*-algebras},
Trans. Amer. Math. Soc. \textbf{132} (1968), 79--99.

\bibitem[Ca]{Ca}
J. W. Calkin,
\emph{Two-sided ideals and congruences in the ring of bounded operators in Hilbert space},
Ann. of Math. (2) \textbf{42} (1941), 839--873.

\bibitem[C]{C1973}
A. Connes,
\emph{Une classification des facteurs de type III},
Ann. Sci. \'{E}cole Norm. Sup. (4) \textbf{6} (1973), 133--252.

%\bibitem[C2]{C1976}
%A. Connes,
%\emph{Classification of injective factors: cases $\text{II}_1$, $\text{II}_\infty$, $\text{III}_\lambda$, $\lambda \neq 1$},
%Ann. of Math. (2) \textbf{104} (1976), 73--115.

\bibitem[D1]{D}
J. Dixmier,
\emph{Sous-anneaux ab\'{e}liens maximaux dans les facteurs de type fini},
Ann. of Math. (2) \textbf{59} (1954), 279--286.

\bibitem[D2]{Dbook}
J. Dixmier,
Les Alg\`{e}bres d'Op\'{e}rateurs dans l'Espace Hilbertien (Alg\`{e}bres de von Neumann),
Gauthier-Villars, Paris, 1969.

\bibitem[Fa]{Fa}
I. Farah,
\emph{All automorphisms of the Calkin algebra are inner},
arXiv preprint math.OA/0705.3085.

\bibitem[FS]{FS}
W. Feit and G. M. Seitz,
\emph{On finite rational groups and related topics},
Illinois J. Math. \textbf{33} (1989), 103--131.

\bibitem[F]{F}
P. A. Fillmore,
\emph{Extensions relative to semifinite factors},
in: Symposia Mathematica, Vol. XX (Convegno sulle Algebre $C^*$ e loro Applicazioni in Fisica Teorica, Convegno sulla Teoria degli Operatori Indice e Teoria $K$, INDAM, Roma, 1975), pp. 487--496, Academic Press, London, 1976.

\bibitem[GN]{GN}
I. Gelfand and M. Neumark,
\emph{On the imbedding of normed rings into the ring of operators in Hilbert space},
Rec. Math. [Mat. Sbornik] N.S. \textbf{12 (54)} (1943), 197--213.

\bibitem[G]{G}
P. Golberg,
\emph{The Sylow $\Pi$-subgroups of locally normal groups} (Russian),
Rec. Math. [Mat. Sbornik] N.S. \textbf{19 (61)} (1946), 451--460.

\bibitem[HS1]{HSeq}
U. Haagerup and E. St\o rmer,
\emph{Equivalence of normal states on von Neumann algebras and the flow of weights},
Adv. Math. \textbf{83} (1990), 180--262.

\bibitem[HS2]{HSpinn}
U. Haagerup and E. St\o rmer,
\emph{Pointwise inner automorphisms of von Neumann algebras},
J. Funct. Anal. \textbf{92} (1990), 177--201.

\bibitem[HS3]{HSaut}
U. Haagerup and E. St\o rmer,
\emph{Automorphisms which preserve unitary equivalence classes of normal states}, in: Operator Theory: Operator Algebras and Applications (Durham, 1988), pp. 531--537, Proc. Sympos. Pure Math., vol. 51, pt. 1, Amer. Math. Soc., Providence, 1990.

\bibitem[HS4]{HSinj}
U. Haagerup and E. St\o rmer,
\emph{Pointwise inner automorphisms of injective factors},
J. Funct. Anal. \textbf{122} (1994), 307--314.

%\bibitem[H]{H1981}
%D. W. Hadwin,
%\emph{Nonseparable approximate equivalence},
%Trans. Amer. Math. Soc. \textbf{266} (1981), 203--231.

\bibitem[IPP]{IPP}
A. Ioana, J. Peterson, and S. Popa,
\emph{Amalgamated free products of $w$-rigid factors and calculation of their symmetry groups},
Acta Math., to appear.

\bibitem[Je]{Je}
T. Jech,
Set Theory (3rd ed.),
Springer-Verlag, Berlin, 2003.

\bibitem[J]{J}
B. E. Johnson,
\emph{Local derivations on $C^*$-algebras are derivations},
Trans. Amer. Math. Soc. \textbf{353} (2001), 313--325.

\bibitem[K1]{K1966}
R. V. Kadison,
\emph{Derivations of operator algebras},
Ann. of Math. \textbf{83} (1966), 280--293.

\bibitem[K2]{K1990}
R. V. Kadison,
\emph{Local derivations},
J. Algebra \textbf{130} (1990), 494--509.

\bibitem[L]{L}
E. C. Lance,
\emph{Automorphisms of certain operator algebras},
Amer. J. Math. \textbf{91} (1969), 160--174.

\bibitem[LS]{LS}
D. R. Larson and A. R. Sourour,
\emph{Local derivations and local automorphisms of $\B(X)$},
in Operator Theory: Operator Algebras and Applications (Durham, 1988), pp. 187--194, Proc. Sympos. Pure Math., vol. 51, pt. 2, Amer. Math. Soc., Providence, 1990.

\bibitem[N]{N}
B. H. Neumann,
\emph{Not quite inner automorphisms},
Bull. Austral. Math. Soc. \textbf{23} (1981), 461--469.

\bibitem[Pe]{Pe}
C. Pearcy,
\emph{$W^*$-algebras with a single generator},
Proc. Amer. Math. Soc. \textbf{13} (1962), 831--832.

\bibitem[PW]{PW}
N. C. Phillips and N. Weaver,
\emph{The Calkin algebra has outer automorphisms},
Duke Math. J. \textbf{139} (2007), 185--202.

\bibitem[P1]{Pkad}
S. Popa,
\emph{On a problem of R. V. Kadison on maximal abelian *-subalgebras in factors},
Invent. Math. \textbf{65} (1981), 269--281.

\bibitem[P2]{Porth}
S. Popa,
\emph{Orthogonal pairs of *-subalgebras in finite von Neumann algebras},
J. Operator Theory \textbf{9} (1983), 253--268.

\bibitem[P3]{Psing}
S. Popa,
\emph{Singular maximal abelian *-subalgebras in continuous von Neumann algebras},
J. Funct. Anal. \textbf{50} (1983), 151--166.

%\bibitem[Sah]{Sah}
%C.-H. Sah,
%\emph{Automorphisms of finite groups},
%J. Algebra \textbf{10} (1968), 47--68.

\bibitem[Sa]{Sa}
S. Sakai,
\emph{Derivations of $W^*$-algebras},
Ann. of Math. \textbf{83} (1966), 273--279.

\bibitem[\v{S}e]{Se}
P. \v{S}emrl,
\emph{Local automorphisms and derivations on $B(H)$},
Proc. Amer. Math. Soc. \textbf{125} (1997), 2677--2680.

%\bibitem[S1]{S2005}
%D. Sherman,
%\emph{Unitary orbits of normal operators in von Neumann algebras},
%J. Reine Angew. Math. \textbf{605} (2007), 95--132.

\bibitem[S]{S2006div}
D. Sherman,
\emph{Divisible operators in von Neumann algebras},
arXiv preprint math.OA/0611364.

\bibitem[T]{T1973}
M. Takesaki,
\emph{Duality for crossed products and the structure of von Neumann algebras of type III},
Acta Math. \textbf{131} (1973), 249--310.

\bibitem[Wa]{Wa}
G. E. Wall,
\emph{Finite groups with class-preserving outer automorphisms},
J. London Math. Soc. \textbf{22} (1947), 315--320 (1948).

\bibitem[W]{W}
W. Wogen,
\emph{On generators for von Neumann algebras},
Bull. Amer. Math. Soc. \textbf{75} (1969), 95--99.

\end{thebibliography}
\end{document}